\documentclass[12pt]{article}

\usepackage{amsfonts}

\newtheorem{ prop}{Proposition}[section]

\newtheorem{thm}{Th\'eor\`eme}[section]

\newtheorem{prop}{Proposition}[section]

\begin{document}

\title{R\'ealisations de Big\`ebres (II) :\\ Th\'eor\`eme de Dualit\'e.}
\author{Eric Mourre}

\maketitle

\bigskip\bigskip

\begin{abstract}
In [1] a method for realizations of bialgebras is presented . It allows
to construct algebras of linear operators wich are endowed with a bialgebra structure :
 $U_x (\Delta ,\epsilon )$ . The inital data for such a construction is a linear 
application :
$$ x:L \rightarrow Invd(F) \ , $$
where $L$ and $F$ are two coalgebras that are supposed of finite dimensions in this
article , and where
$Invd(F) $  denotes the algebra of right invariant operators on any coalgebra $F$ .
 $T(L) $ and $T(F) $ denote the tensor algebras which are equiped with their bialgebra
structures extending, by algebra morphisms, the coproducts and counits defined on $L$,
resp on
$F$ . Let us recall shortly the construction . \\
The linear application $x :L \rightarrow Invd(F) $ and the coalgebra structure of $L$
allow to extend for any $l\in L$ , $x(l)$ into an operator $X(l)$ not only defined on
$F$ , but defined on $T(F)$ in the following manner :\\ for $w_1 $ and $w_2 $ in $T(F) $
,
$$ X(l) (w_1 . w_2 ) = \sum_{k} X(l_{k}^{'})(w_1) . X(l_{k}^{''})(w_2) $$
where
$$\Delta_L (l) = \sum_k l_k^{'}\otimes l_k^{''} \ .$$
In [1] it is shown that the algebra generated in $ End(T(F))$ by the operators $X(L)$ and
 the indentity 
is a bialgebra : $U_x (\Delta_L , \epsilon_L )$ , which is included in $Invd(T(F))$ .\\
Then the duality theorem is expressed in the following way :\\
By taking the transposition of the application :
$$x :L\rightarrow Invd(F)  $$
we obtain :
$$x^t = y :F\rightarrow Invd(L) \ ,$$
and by applying the previous results to $ x$ and $y$  we obtain two bialgebras $U_x
(\Delta_L,\epsilon_L)
$ and
$V_y (\Delta_F ,\epsilon_F )$ , and two bialgebra morphisms mapping units onto units
$$\pi_x :T(L) \rightarrow U_x \subset Invd(T(F)) \subset End(T(F)) \ ,$$
 
$$\pi_y : T(F) \rightarrow V_y \subset Invd(T(L)) \subset End(T(L)) \ .$$
 By denoting for any coalgebra $C$,
the canonical duality between  $c\in C$ and
$A\in Invd(C)$ : $\epsilon_C \circ A (c)  = <A,c>$ , then we have :
\begin{thm}{ Duality theorem }\\

 For
any
$w\in T(L) $ and $v\in T(F)$ the following indentity holds :
$$ <\pi_x (w) ,v> = < \pi_y (\tau (v)),\tau (w)> \ ,$$
where $\tau $ denotes the anti-isomorphism of tensor algebras which reverses the order
of tensors  . 
\end{thm}
Then this theorem is used to prove properties of the ideal of relations associated to
any realized bialgebra $U_x(\Delta_L ,\epsilon_L )$ . \\In particular the vector spaces
of relations of order less or equal to $n$ in $T(L)  $, are shown to be coideals , for
any
$n$ . Then we describe two complementary constructions of these coideals . \\
One of these constructions relies on the following characterization of the
vector space of relations of $U_x $ contained in a subcoalgebra $C$ of $T(L)$ : it
coincides with the maximal coideal of $C$ contained  in the set of relations defined 
by the minimal representation, on $\mathbb{C} \oplus F$ .

\end{abstract}
\vskip0.3cm
CPT 2005 /P.004 \\[0.1cm]
Centre de Physique Th\'eorique , CNRS , Luminy, case 907 ,
13288 Marseille Cedex 9 , France .

\section{Introduction .}

L'article " Remarques sur une m\'ethode de r\'ealisations de big\`ebres, et alg\`ebres
de Hopf associ\'ees \`a certaines r\'ealisations " [1] , pr\'esente une m\'ethode
g\'en\'erale, qui permet de construire des alg\`ebres d'op\'erateurs lin\'eaires, qui
sont munies d'une structure de big\`ebre.\\
Une big\`ebre est une alg\`ebre associative (sur $\mathbb{C}$), $U$ , munie d'un
coproduit coassociatif ,\\
$\Delta :U
\rightarrow U\otimes U$ , qui est un morphisme d'alg\`ebres , et d'une counit\'e\\
$\epsilon : U\rightarrow \mathbb{ C}$ , qui est un morphisme d'alg\`ebres , et qui
v\'erifient :\\
$$\epsilon \otimes id \circ \Delta =id\otimes \epsilon  \circ \Delta = id : U\rightarrow
U
\ .$$  Soient deux cog\`ebres , $L(\Delta_L ,\epsilon_L ) $ et $ F(\Delta_F ,
\epsilon_F)$ ,  les alg\`ebres tensorielles $T(L)$ et $T(F)$ sont munies de structures
de big\`ebres qui
\'etendent par morphismes d'alg\`ebres les structures de cog\`ebres de $L$
et respectivement de $F$ .\\
Dans le cas g\'en\'eral de cog\`ebres $L$ et $F$ de dimensions quelconques ,voir(1), il
y a des pr\'ecautions \`a  prendre pour d\'efinir des sous alg\`ebres d'op\'erateurs 
invariants \`a droite sur les cog\`ebres  $L$ , $F$ , et surtout sur, $T(L)$ et $T(F)$
.\\
\\ 
Dans cet article on suppose que les cog\`ebres $L$ et $F$ sont de dimensions finies et
pour fixer les notations on notera respectivement  $K$ et $E$ les
alg\`ebres associatives , avec unit\'es , dont $L$ , respectivement $F$ , sont les
cog\`ebres duales:
$$L(\Delta_L ,\epsilon_L )=K^* \ \ et\ \  F(\Delta_F ,\epsilon_F) =E^* \ .$$
\\
Alors les alg\`ebres d'op\'erateurs invariants \`a droite sont bien d\'efinies ;
on les d\'esignera par $Invd(L)$, $Invd(F)$ , $Invd(T(L))$ , et $Invd(T(F))$ .
Ces alg\`ebres sont vectoriellement isomorphes aux duaux des cog\`ebres
correspondantes.\\
\\
Les donn\'ees initiales pour la r\'ealisation d'une big\`ebre sont :\\
\\
deux cog\`ebres  $L(\Delta_L ,\epsilon_L ) $ et $ F(\Delta_F , \epsilon_F)$ et une
application lin\'eaire $x$ :
$$x :L \rightarrow Invd(F)\ . $$
Alors cette application lin\'eaire $x$ et la structure de cog\`ebre de $L(\Delta_L ,
\epsilon_L)$ \\ permettent, pour tout vecteur $l\in L$ , d' \'etendre de mani\`ere
unique l'op\'erateur $x(l)$ en un op\'erateur $X(l)$ , non plus simplement de :
 $F\rightarrow F$ ,  mais \\    de : $T(F) \rightarrow T(F)$ , de la fa\c{c}on suivante
:\\ pour $w_1$ , $w_2$ dans $T(F)$  :
$$X(l)(w_1 .w_2 ) = \sum_k X(l_{k}^{'})(w_1).X(l_{k}^{''})(w_2) $$
o\`u 
$$\Delta (l) = \sum_k l_{k}^{'}\otimes l_{k}^{''} \ . $$
L'alg\`ebre engendr\'ee par les op\'erateurs  $X(l) ,\ l\in L$ et l'op\'erateur
identit\'e  dans
$End(T(F))$ est alors munie d'une structure de big\`ebre not\'ee:\\
$$U_x(\Delta_L ,\epsilon_L)\ ;\ avec \ : \  \Delta_L X(l) =\sum_k X(l_k^{'})\otimes
X(l_k^{''}) \ ;$$  elle agit sur $T(F)$.\\
\\
Les op\'erateurs $X(l) \in Invd(T(F))$ sont des op\'erateurs invariants \`a droite sur
$T(F)$ :
$$U_x(\Delta_L ,\epsilon_L)\subset Invd(T(F)) \subset End(T(F)).$$
De plus l'on obtient alors un morphisme d'alg\`ebres pr\'eservant l'unit\'e :
$$\pi_{x} :T(L) \rightarrow U_x(\Delta_L ,\epsilon_L) \subset Invd(T(F)) $$
qui est aussi un morphisme de cog\`ebres.\\
\\
Le paragraphe 2 rappelle les th\'eor\`emes d\'ecrivant la r\'ealisation de big\`ebres 
de [1] dans le cas o\`u $L$ et $ F$ sont de dimensions finies .\\
\\
Le paragraphe 3 est consacr\'e \`a la d\'emonstration du th\'eor\`eme de dualit\'e
que l'on d\'ecrit maintenant . \\
\\
La construction pr\'ec\'edente est associ\'ee \`a toute application lin\'eaire\\
 $x$ d'une cog\`ebre $L$ dans l'espace des op\'erateurs invariants \`a droite sur une
cog\`ebre $F$ .\\ Parce que le dual d'une cog\`ebre de dimension finie est canoniquement
isomorphe \`a l'espace vectoriel des op\'erateurs invariants \`a droite sur cette 
cog\`ebre , par transposition de l'application $x:L\rightarrow Invd(F) $ on obtient, 
l'application lin\'eaire $ y:F\rightarrow Invd(L)$ : 
 $$ x^t = y :F\rightarrow Invd(L) \ . $$  
On obtient ainsi :\\
Une application lin\'eaire 
$$Y :F\rightarrow Invd(T(L)) \ ,$$ 
la big\`ebre :
$$V_y (\Delta_{F},\epsilon_{F}) \subset Invd(T(L)) $$
et un morphisme d'alg\`ebres $\pi_y$ :
$$\pi_y : T(F) \rightarrow V_y(\Delta_{F},\epsilon_F ) \subset Invd(T(L)) \ ,$$
qui est aussi un morphisme de cog\`ebres .\\
\\
Pour tout op\'erateur $a \in Invd(C) $ et tout $c\in C$ , o\`u $C$ est une cog\`ebre on
d\'enote la dualit\'e entre $a$ et $c $ par : $<a,c>$ ,

$$ \epsilon_{C} \circ a (c) = <a,c> \ .$$
En d\'esignant , pour tout espace vectoriel $V$ par $T_0(V)$ la sous alg\`ebre
tensorielle suivante ,
$$T_0 (V) = V\oplus V\otimes V \oplus ...\oplus \otimes^{n} V \oplus ...\ , $$
et par $\tau $ l'anti-isomorphisme d'alg\`ebre tensorielle correspondant au renversement
de l'ordre des tenseurs ,\\
\\
alors le th\'eor\`eme de dualit\'e s'exprime de la fa\c{c}on suivante :
\begin{thm}{Th\'eor\`eme de Dualit\'e.}\\
Pour toutes suites $(l_{1},l_{2},...l_{p}) $, $l_{i}\in L$ et  $(f_{1},f_{2},...f_{q})\ ,
$  $f_{j}\in F$, nous avons :\\
$$< X(l_{1})\circ X(l_{2})\circ...\circ X(l_{p}) , f_{1} \otimes f_{2}\otimes ...\otimes
f_{q}) >$$ est \'egal \`a

$$< Y(f_{q})\circ Y(f_{q-1})\circ...\circ Y(f_{1}) , l_{p}\otimes
l_{p-1}\otimes...\otimes l_{1}>\ .$$ 
Et plus g\'en\'eralement , pour tous vecteurs $w\in T(L)$ et $v\in T(F)$
l'on a :
$$< \pi_x (w),v> = <\pi_y (\tau (v)), \tau(w) > \ .$$

\end{thm}
\vskip 0.5 cm
Le paragraphe (4) est consacr\'e \`a l'\'etude de propri\'et\'es concernant les id\'eaux
de  relations des big\`ebres r\'ealis\'ees .
En particulier il y est d\'emontr\'e que pour tout $n$ , l'espace vectoriel, $R_n
(U_x)$ , des relations d'ordre $n$ d'une big\`ebre $U_x(\Delta_L ,\epsilon_L)$ est un
coid\'eal de
$T(L)$ .\\
\\
L'espace vectoriel $R_n(U_x)$ des relations d'ordre $n$ de la big\`ebre
$U_x(\Delta_L,\epsilon_L)$
\'etant le noyau de l'application:
$$ \pi_x : T_n(L) \rightarrow U_x \subset Invd(T(F)) \ ,$$
o\`u 
$$T_n(L): \mathbb{C} \oplus L\oplus L\otimes L\oplus  ..\oplus \otimes^n L \ .$$
La propri\'et\'e de $R_n(U_x)$ d'\^etre un coid\'eal de $T_n (L)$ est d\'emontr\'ee dans
le th\'eor\`eme 4.1 :
\begin{thm}
Pour toute big\`ebre r\'ealis\'ee , $U_x(\Delta_L ,\epsilon_L)$ , et pour tout  $n$ ,
l'on a :
$$\Delta_L (R_n(U_x)) \subset R_n(U_x)\otimes T_n(L) + T_n(L) \otimes R_n(U_x) \ . $$

\end{thm}
D'autre part on obtient , toujours comme cons\'equences du th\'eor\`eme de dualit\'e,
deux constructions , compl\'ementaires et utiles , des coid\'eaux de relations d'ordre
$n$ .\\
\\
La premi\`ere construction provient de la caract\'erisation des coid\'eaux de relations
de la big\`ebre
$U_x(\Delta_L ,\epsilon_L) $ comme les coid\'eaux , de relations d\'efinies par la
repr\'esentation minimale , en l'occurrence la repr\'esentation de $U_x$ 
sur $ \mathbb{C} \oplus F$ ( th\'eor\`emes 4.3 , 4.4 , paragaphe 4.2 ) . \\
Cette caract\'erisation des coid\'eaux de relations de la big\`ebre $U_x(\Delta_L
,\epsilon_L)$ est int\'eressante ; par exemple elle pourrait permettre de trouver des
solutions au 
probl\`eme de  la r\'ealisation de big\`ebres admettant l'action d'op\'erateurs 
particuliers .\\ Pour la position du probl\`eme de
l'existence de l'antipode sur une big\`ebre $U_x(\Delta_L ,\epsilon_L )$ , lorsque  par
exemple  
$K$ est trigonalisable , voir [1] . \\
  \\ Cette premi\`ere construction des coid\'eaux
$R_n(U_x)$ est d\'emontr\'ee dans\\ le th\'eor\`eme 4.4 et consiste dans les
op\'erations suivantes .\\
\\
Soit
le sous espace vectoriel $R_n^{1} $ de $T_n (L)$ d\'efini par le noyau de $\pi_x^{1}$ :
$$ R_n^{1} = T_n (L)\cap ker\  \pi_x^{1}\ \ ,\ \pi_x^{1} :T(L)\rightarrow
Invd(F)\subset Invd(T(F)) \ .$$ 
L'alg\`ebre  $T_n (K)$ est le dual de $T_n (L)$ et l'orthogonal dans $T_n(K)$ de
$R_n^{1}$ est un sous espace vectoriel $W_n \subset T_n (K)$ de l'alg\`ebre $T_n (K)$ .\\
Soit $ B(W_n)$ la plus petite sous alg\`ebre de $T_n (K)$ contenant $W_n$ et l'unit\'e
de $T_n (K) $;\\ on a alors :
\begin{thm}
L'orthogonal dans $T_n (L)$ de la sous alg\`ebre $B(W_n)$ est $R_n(U_x)$, le coid\'eal
des relations d'ordre $n$ . 

\end{thm}
Finalement , dans cet article , on expose une construction compl\'ementaire, issue du
th\'eor\`eme 4.1 , des coid\'eaux $R_n(U_x) $ dans le th\'eor\`eme 4.5 .\\
\\
De plus l'on 
remarque dans ce paragraphe 4 , que les r\'esultats restent valables, non plus
simplement pour les sous espaces vectoriels des relations de la big\`ebre
$U_x(\Delta_L ,\epsilon_L)$ contenues dans les cog\`ebres $T_n(L)$  , mais pour
 l' espace vectoriel des relations contenues dans une sous cog\`ebre(de dimension
finie),  arbitraire , de $T(L)$ .\\
\\
Pour les propri\'et\'es et les d\'efinitions usuelles des alg\`ebres, cog\`ebres, et 
des alg\`ebres de Hopf voir [2] . En ce qui concerne les constructions des groupes
quantiques et leurs \'etudes , ainsi que les applications et leurs origines voir
[3],[4],[5] et les r\'ef\'erences qui y sont contenues , dont l'origine en termes
de  $ C^{*}$ Alg\`ebres. D'autre part des travaux concernent l'utilisation dans le
domaine de l'analyse fonctionnelle et de la g\'eom\'etrie non commutative d'alg\`ebres
de Hopf particuli\`eres [6], ou encore dans le domaine de la th\'eorie quantique des
champs [7] , [8]  . Signalons enfin une approche, plus axiomatique , 
des alg\`ebres de Hopf en termes d'alg\`ebres de Von Neumann.

\section{Rappels sur une m\'ethode de r\'ealisations de
big\`ebres.}
 Dans ce paragraphe on d\'ecrit la m\'ethode de r\'ealisation de big\`ebres
plubli\'ee dans: arXiv:math-ph 0312049.\\ Les donn\'ees
initiales pour la r\'ealisation d'une big\`ebre sont les suivantes.\\
Soient $K$ et $E$ deux alg\`ebres associatives avec unit\'es ;\\
on d\'esigne par $L$ et $F$ les cog\`ebres (coassociatives avec counit\'es) duales de
$K$ et$E$. Pour tout espace vectoriel $F$, $T(F)$ d\'esigne l'alg\`ebre tensorielle
construite sur $F$. Dans cet article on supposera pour simplifier que les alg\`ebres \\
$K$ et $E$ sont de dimensions finies sur $\mathbb{C}$ .\\
Alors pour toute application lin\'eaire :$$x_{t} : L\rightarrow E $$
est associ\'ee une big\`ebre $ U_{x}(\Delta_{L},\epsilon_{L})$, construite de la
fa\c{c}on suivante :\\
pour tout $l\in L$ , $x_{t}(l) \in E$ ; et la transposition de l'op\'erateur
multiplication \`a gauche dans $E$ par $x_{t}(l) $ est un op\'erateur invariant \`a
droite sur la cog\`ebre $F$. On d\'esigne par $ x$ l'application lin\'eaire de la
cog\`ebre $L$ dans les op\'erateurs invariants \`a droite sur la cog\`ebre $F$ .
L'application lin\'eaire $x$ et la structure de cog\`ebre de $L$ permettent d'\'etendre
les op\'erateurs $x(l)$ non plus simplement de $ F\rightarrow F$ mais de
$T(F)\rightarrow T(F)$ de la mani\`ere suivante :\\
pour $w_{1} $ et $w_{2}$  $\in T(F) $ on d\'efinit :\\

$$ X(l)(w_{1}.w_{2}) = \sum_{k}X(l_{k}^{'})(w_{1}).X(l_{k}^{''})(w_{2}) $$ 
o\`u
$$ \Delta_L (l) = \sum_{k}l_{k}^{'}\otimes l_{k}^{''}$$
On obtient ainsi une application lin\'eaire :

$$ X:L \rightarrow End(T(F))$$
et pour toute base $(l_{i})$ de $L$ , la sous alg\`ebre de $End(T(F))$ engendr\'ee par
les op\'erateurs $X(l_{i})$ et l'identit\'e est alors munie d'une structure de big\`ebre,
que l'on d\'esigne par $U_{x}(\Delta_{L},\epsilon_{L})$ .\\
On rappelle dans ce paragraphe quelques points essentiels de la construction dont nous
aurons besoin dans le paragraphe 3 pour formuler et \'etablir le th\'eor\`eme de
dualit\'e et dans le paragraphe 4 pour d\'emontrer ses applications concernant la
d\'escription de propri\'et\'es de l'id\'eal des relations dans  la big\`ebre 
$U_{x}(\Delta_{L},\epsilon_{L})$. \\

\begin{prop}
 Pour toute cog\`ebre $F$ de dimension finie, de coproduit $\Delta_{F}$, et counit\'e
$\epsilon_{F}$, $T(F)$ est munie d'une structure de big\`ebre not\'ee:
$T(F)(\Delta_{F},\epsilon_{F})$.\\ 
Soit $Invd(T(F))$  la sous alg\`ebre de $ End(T(F))$, form\'ee par les op\'erateurs $X$,
invariants \`a droite, v\'erifiant :\\
$\Delta_{F}\circ X = X \otimes id \circ \Delta_{F} \ .$\\
L'application lin\'eaire $$ X\in Invd(T(F)) \rightarrow \omega_{X}=\epsilon_{F} \circ X
\in T(F)^{*}$$
est un anti-morphisme d'alg\`ebres; l'application inverse est donn\'ee par:
$$X = \omega_{X}\otimes id \circ \Delta_{F}\ \ .$$

\end{prop}
Nous avons :
$$T(F)=C\oplus F\oplus F\otimes F\oplus\ldots\oplus\otimes^{n}F \oplus\ldots$$
Le coproduit $\Delta_{F}$ et la counit\'e $\epsilon_{F}$ sont des morphismes d'alg\`ebres
respectivement de $T(F) \rightarrow T(F)\otimes T(F)$ et de $T(F) \rightarrow \mathbb{C}$
.\\ Pour $1\in \mathbb{C}$ , $\Delta_{F}(1)=1\otimes 1$ et $\epsilon_{F} (1)= 1$;\\ pour
$f
\in F
\subset T(F) $ , le coproduit $\Delta_{F}$, co\"{\i}ncide avec le coproduit sur $F$ ,
ainsi que la counit\'e . \\ 
$F$ \'etant le dual d'une alg\`ebre de dimension finie $E$  par hypoth\`ese, en prenant
une base  
$e_{i}$ pour $E$ on obtient la base canonique duale  $f_{i}$ pour $F$ et donc la base 
canonique associ\'ee sur $T(F)$ .\\
Si on note :$$\Delta f_{i} = \sum_{k}f_{i,k}^{'}\otimes f_{i,k}^{''}$$
on a :
$$ \Delta_{F} (f_{i1}\otimes f_{i2}\otimes\ldots\otimes f_{in}) =
\sum_{k1,k2,...,kn}(f_{i1,k1}'\otimes f_{i2,k2}'\otimes ..\otimes f_{in,kn}')\otimes
(f_{i1,k1}''\otimes  f_{i2,k2}''\otimes ..\otimes f_{in,kn}'')
$$
$\epsilon_{F} :T(F) \rightarrow \mathbb{C} $ \'etant un morphisme d'alg\`ebres de $ T(F)
\rightarrow \mathbb{C}$, on a bien par exemple :
 $$  \epsilon_{F} (f_{i1}\otimes f_{i2}\otimes\ldots\otimes f_{in}) = \epsilon_{F}
(f_{i1}).\epsilon_{F} (f_{i2})\ldots\epsilon_{F} (f_{in})) \  ,$$

$$ id \otimes \epsilon_{F}  (\Delta_{F} (f_{i1}\otimes
f_{i2}\otimes\ldots\otimes f_{in})) =  f_{i1}\otimes f_{i2}\otimes\ldots\otimes f_{in} \
.
$$
$T(F)$ \'etant une big\`ebre l'alg\`ebre des op\'erateurs invariants \`a droite 
$$Invd (T(F)(\Delta_{F},\epsilon_{F})) \subset End(T(F)) $$
est l'alg\`ebre des op\'erateurs $X \in End(T(F)) $ qui v\'erifient :
$$ \Delta_{F} \circ X = X\otimes id \circ \Delta_{F}\ \ .$$ 
Les formes $\omega \in T(F)^{*}$ \'etant munies du produit associatif :
$$\omega_{1}.\omega_{2} = \omega_{1} \otimes \omega_{2} \circ \Delta_{F} \ , $$
alors l'application lin\'eaire :
$$ i : X\in Invd(T(F)) \rightarrow i(X) = \epsilon_{F} \circ X \in T(F)^{*} $$
est un anti-morphisme d'alg\`ebres .\\
Cette dualit\'e  entre les op\'erateurs $X \in Invd(T(F))$ et $ c\in T(F)$ sera souvant
not\'ee par:$$<X,c>= \epsilon_{F}(X(c)) \ \ .$$
Dans les th\'eor\`emes suivants on d\'ecrit des r\'esulats obtenus dans [1].\\

\begin{thm}
Soient deux cog\`ebres $L(\Delta_{L},\epsilon_{L})$ et $F(\Delta_{F},\epsilon_{F})$,
et une application lin\'eaire $x$,
$$x : L\rightarrow Invd(F) $$
alors il existe une unique application lin\'eaire $X$ :
$$ X: L \rightarrow  Invd(T(F))$$
telle que: \\[0.2cm]
1) pour $1 \in T(F)$ $\ ,$ $ X(l)(1) = \epsilon (l) 1 $ , \\[0.2cm]
2) pour $f \in F \subset T(F) $ \ , $ X(l)(f) = x(l)(f) $ \ , \\[0.2cm]
3) pour tout $w_{1} , w_{2} \in T(F) $ alors :
$$ X(l)(w_{1}.w_{2}) = \sum_{k} X(l_{k}^{'})(w_{1}). X(l_{k}^{''})(w_{2}) \ .$$
4) Les op\'erateurs $X(l) $ sont des op\'erateurs invariants \`a droite sur $T(F)$ , 
et ils v\'erifient :\\[0.2cm]
5) $ X(l)( \otimes^{n}F) \subset \otimes^{n}F $

\end{thm}
La d\'emonstration de ce th\'eor\`eme est donn\'ee dans [1 ] dans un cadre 
un peu plus g\'en\'eral que celui pr\'esent\'e ici, d'une application lin\'eaire $x$,
d'une cog\`ebre $L$ de dimension finie, \`a valeurs dans les invariants \`a droite 
sur une cog\`ebre $F$ de dimension finie.\\[0.4cm]

\begin{thm}
Soient deux cog\`ebres $L(\Delta_{L},\epsilon_{L})$ et $F(\Delta_{F},\epsilon_{F})$,
et une application lin\'eaire $x$,
$$x : L\rightarrow Invd(F) $$
et $X$ l'unique application lin\'eaire :
$$ X: L \rightarrow  Invd(T(F))$$
d\'efinie dans le th\'eor\`eme pr\'ec\'edent.\\[0.2cm]
Alors l'alg\`ebre $U_x$ engendr\'ee par les op\'erateurs $X(l)\in Invd(T(F)) $ et
l'identit\'e  est munie d'une structure de big\`ebre : $U_{x}(\Delta_{L},\epsilon_{L})$ ,
le coproduit et la counit\'e qui sont des morphismes d'alg\`ebres sont d\'efinis sur les
g\'en\'erateurs par:\\
$$\Delta_{L}(X(l)) = X\otimes X \circ \Delta_{L}(l) \ \ .$$ 
$$\epsilon_{L}(X(l))=\epsilon_{L}(l) $$  
 
\end{thm}
La d\'emonstration est donn\'ee dans [1] ; le point essentiel qui permet de montrer 
que le coproduit $ \Delta_{L}: U_{x} \rightarrow     U_{x}\otimes  U_{x}$ est bien
d\'efini, repose sur le fait que les op\'erateurs $X(l)$ sont des op\'erateurs
invariants \`a droite sur la big\`ebre $T(F)$ . 
En ce qui concerne la counit\'e de $U_{x}$ on a par exemple, pour tout op\'erateurs $u\in
U_{x}$ :
$$ \epsilon_{L}(u).1= u(1)\ ;\ \ et \ \ u(1. w)= \sum_{k} u_{k}^{'}(1).u_{k}^{''}(w)\ \ 
.$$
pour tout $w \in T(F))$, ce qui d\'emontre bien:  $\epsilon_{L}\otimes id \circ
\Delta_{L}(u) = u $ .

\section{ Th\'eor\`eme de Dualit\'e .}

Soient $L(\Delta_{L},\epsilon_{L})$ et $F(\Delta_{F},\epsilon_{F})$ 
des cog\`ebres de dimensions finies, avec counit\'es. Ces deux cog\`ebres sont les
duaux des alg\`ebres associatives $K$ et$E$ .\\
La donn\'ee de l'application lin\'eaire :
$$x :L\rightarrow Invd(F)$$ 
est \'equivalente \`a la donn\'ee d'une application lin\'eaire $x_{t}$ :

$$x_{t} : L \rightarrow E $$
et la transposition de cette application donne une application $y_{t}$ :

$$ y_{t} :F \rightarrow K $$ et donc une application lin\'eaire $y$ :

$$ y: F \rightarrow Invd(L)\ \ . $$
Ainsi en appliquant les th\'eor\`emes pr\'ec\'edents aux applications
$x$ et $y$ on obtient les deux applications lin\'eaires $X$ et $Y$ :

$$X :L \rightarrow Invd(T(F)) $$
$$Y :F \rightarrow Invd(T(L)) $$
 les big\`ebres :
$$U_{x}(\Delta_{L},\epsilon_{L}) \subset Invd(T(F))      $$ 
$$V_{y}(\Delta_{F},\epsilon_{F}) \subset Invd(T(L)) \ \ .    $$ 
et deux morphismes de big\`ebres  $\pi_x $ et $\pi_y$ pr\'eservant l'unit\'e ,
$$\pi_x :T(L)\rightarrow U_x(\Delta_L ,\epsilon_L ) \subset Invd(T(F))\subset End(T(F))$$
$$\pi_y :T(F) \rightarrow V_y(\Delta_F ,\epsilon_F) \subset Invd(T(L)) \subset
End(T(L)) \ .$$ Le th\'eor\`eme suivant met en \'evidence une relation de dualit\'e entre
ces deux big\`ebres .\\
\\
{\bf{D\'efinitions}}\\
Pour tout \'el\'ement $c$ d'une cog\`ebre $C$ et tout op\'erateur $A$
invariant \`a droite sur la cog\`ebre $C$ , on note :
$$<A,c> = \epsilon_{C} \circ A (c)  \ \ .$$
De plus, pr\'ecisons ici que dans le th\'eor\`eme suivant, les cog\`ebres sont
respectivement $T(F)$ et $T(L)$ et que les op\'erateurs invariants \`a droite
concern\'es sont ceux qui sont  respectivement dans $U_{x}$ et $V_{y}$ .\\
 \\
Pour toute alg\`ebre tensorielle sur $V$ on d\'esigne par $\tau $ l'anti-isomorphisme
d'alg\`ebre tensorielle correspondant au renversement de l'ordre des tenseurs, et \\
$$T_0 (V) = V\oplus V\otimes V \oplus ...\oplus \otimes^{n}V \oplus ...\ \  .$$

\begin{thm}

Pour toutes suites $(l_{1},l_{2},...l_{p}) $, $l_{i}\in L$ et  $(f_{1},f_{2},...f_{q})\ ,
$\\ $f_{j}\in F$, nous avons :
$$< X(l_{1})\circ X(l_{2})\circ...\circ X(l_{p}) , f_{1} \otimes f_{2}\otimes ...\otimes
f_{q}) >$$ est \'egal \`a

$$< Y(f_{q})\circ Y(f_{q-1})\circ...\circ Y(f_{1}) , l_{p}\otimes
l_{p-1}\otimes...\otimes l_{1}> \ .$$ 
\\
Et pour $w\in T_0 (L) $ et $v\in T_0 (F)$ l'on a :
$$<\pi_{x}(w) ,v> = <\pi_{y}(\tau (v)),\tau (w)> \ .$$
De plus cette identit\'e reste vraie pour tout $w\in T(L)$ et tout $v\in T(F) $ .
\end{thm}
\vskip0.2cm
D\'emonstration.\\
Pour toute application $ x :L\rightarrow Invd(F)$ il existe un unique \'el\'ement
$ w\in K\otimes E$ tel que:
$$ x_{t}(l) =l\otimes id ( w )\ ,$$
et
$$ y_{t}(f) = id\otimes f(w) \ .$$

Nous avons ;\\
$$< X(l),f> =\epsilon_{F}(X(l)(f)) = \epsilon_{F}\circ x(l)\otimes
\epsilon_{F}\circ\Delta_F (f))$$
ce qui s'\'ecrit :
$$<X(l),f>=f(x_t (l)) =l\otimes f(w) $$
De m\^eme nous avons:
$$ <Y(f),l> = \epsilon_{L}(Y(f)(l)) =\epsilon_{L}\circ y(f)\otimes
\epsilon_{L}\circ
 \Delta_{L}(l))= l\otimes f(w) \ . $$
Nous avons donc :

$$<X(l),f> = < Y(f), l>\ . $$
Calculons maintenant:

$$ < X(l_{1})\circ X(l_{2})\circ...\circ X(l_{p}) , f> \ ,$$
ce qui s'\'ecrit :
$$ < X(l_{p})\otimes X(l_{p-1})\otimes...\otimes X(l_{1}) ,\Delta_{F}^{p-1} (f)> \ ,$$
parceque d'apr\'es la proposition (2,1) on sait que l'on a un anti-isomorphisme entre les
op\'erateurs invariants \`a droite sur une cog\`ebre et l'ag\`ebre des formes sur
cette cog\`ebre . D'autre part:
$$<Y(f),  l_{p}\otimes l_{p-1}\otimes...\otimes l_{1}>\ ,$$
est \'egal \`a :
$$ <\otimes^{p}Y \circ \Delta_{F}^{p-1}(f) , l_{p}\otimes l_{p-1}\otimes...\otimes
l_{1}>\ ,$$ 
par construction des op\'erateurs $Y(f)$ et parce que $\epsilon_{L}$ est un
morphisme d'alg\`ebres de $T(L)$ dans $\mathbb{C}$ . Et en utilisant les identit\'es\\
$<X(l),f> = < Y(f), l>$ , l'on obtient :
$$<Y(f),  l_{p}\otimes l_{p-1}\otimes...\otimes l_{1}>=< X(l_{1})\circ
X(l_{2})\circ...\circ X(l_{p}) , f>\ \ . $$
Calculons maintenant :
$$< X(l_{1})\circ X(l_{2})\circ...\circ X(l_{p}) , f_{1}\otimes f_{2}\otimes...\otimes  
f_{q}>\ \ . $$
qui est \'egal \`a 
$$<\Delta_{U_x}^{q-1} (X(l_{1})\circ X(l_{2})\circ...\circ X(l_{p}))\ ,f_{1}\otimes
f_{2}\otimes...\otimes   f_{q}>\ \ . $$
En utilisant les faits que $\Delta_{U_x}^{q-1}$ est un morphisme d'alg\`ebres, ainsi 
que
$\epsilon_{F}$ est un morphisme d'alg\`ebres ($T(F)\rightarrow \mathbb{C} $), et la
formule de d\'efinition des op\'erateurs
$X(l)$ on obtient:

$$<X(l_{1}^{(1)})\circ X(l_{2}^{(1)})..\circ X(l_{p}^{(1)}),f_{1}><X(l_{1}^{(2)})\circ 
X(l_{2}^{(2)})..\circ X(l_{p}^{(2)}),f_{2}>..<X(l_{1}^{q})\circ X(l_{2}^{q})..\circ
X(l_{p}^{q}),f_{q}> $$
o\`u pour chaque $i\in (1,...,p) $ , 
$$\Delta_{L}^{q-1}(l_{i})=l_{i}^{(1)}\otimes l_{i}^{(2)}\otimes...\otimes l_{i}^{(q)}$$
ce qui invoque une sommation bien d\'efinie.\\
\\
 D'autre part calculons:
$$<Y(f_{q})\circ Y(f_{q-1})\circ...\circ Y(f_{1})\ ,\ l_{p} \otimes l_{p-1} \otimes ...
\otimes l_{1}> \ . $$
Les op\'erateurs $Y(f)$ \'etant des op\'erateurs invariants \`a droite sur $T(L)$, on
obtient en utilisant l'anti-isomorphisme entre l'alg\`ebre des op\'erateurs invariants
\`a droite sur $T(L)$ et l'alg\`ebre des formes associ\'ees d\'efinies sur $T(L)
(\Delta_L,\epsilon_L)$:
$$ <Y(f_{1}), l_{p}^{1} \otimes l_{p-1}^{1} \otimes ..
\otimes l_{1}^{1}> <Y(f_{2}), l_{p}^{2} \otimes l_{p-1}^{2} \otimes ..
\otimes l_{1}^{2}>.. <Y(f_{q}), l_{p}^{q} \otimes l_{p-1}^{q} \otimes ..
\otimes l_{1}^{q}>\ \ ,  $$
 en ayant aussi utiliser le fait que $\Delta_{L}^{q-1} $ est un morphisme
d'alg\`ebres .\\ Et en utilisant les identit\'es d\'ej\`a d\'emontr\'ees :

$$<Y(f),  l_{p}\otimes l_{p-1}\otimes...\otimes l_{1}>=< X(l_{1})\circ
X(l_{2})\circ...\circ X(l_{p}) , f>\ \ , $$
on obtient le th\'eor\`eme de dualit\'e , l'extension sur $T(L)$ et $T(F)$ est une
cons\'equence des d\'efinitions .

\section{Dualit\'e et d\'escription de propri\'et\'es des id\'eaux de relations des 
big\`ebres r\'ealis\'ees .}

Dans ce paragraphe il est naturel d'utiliser les notations suivantes.\\
\\
Soient $T(L)$ ou $T(F)$ :
$$T_{0}(L)= L\oplus L\otimes L\oplus\ldots\oplus\otimes^{n}L \oplus\ldots \ ;$$
c'est une sous alg\`ebre et une sous cog\`ebre de $T(L)$ ;
$$T_{n}(L) =\mathbb{C}\oplus  L\oplus L\otimes L\oplus\ldots\oplus\otimes^{n}L           
$$ D'apr\'es les th\'eor\`emes du paragraphe (2) pour tout $w\in K\otimes E $ on 
a un morphisme  $$ \pi_{x} : T(L) \rightarrow U_{x}\subset Invd(T(F)) \subset
End(T(F))$$
donn\'e sur les g\'en\'erateurs par $\pi_{x}(l)= X(l) $et $ \pi_x (1)=id $.\\
Le morphisme d'alg\`ebres $\pi_{x}$ est aussi un morphisme de cog\`ebres:\\
$\pi_{x} :T(L) \rightarrow U_{x}(\Delta_{L},\epsilon_{L})$ .\\
De plus $\pi_{x}$ est non seulement une repr\'esentation de la big\`ebre $T(L)$
sur $T(F)$ mais c'est une action sur cette alg\`ebre:\\
$$ \pi_{x}(z)(w_{1}.w_{2}) = \pi_{x}(z^{1})(w_{1}).\pi_{x}(z^{2})(w_{2}) $$
pour tout $z\in T(L)$ et tout $w_{1}\ ,\ w_{2}\ \in T(F) $ .\\
De m\^eme on a une action de la big\`ebre $T(F)$ sur l'alg\`ebre $T(L)$ , donn\'ee
par le morphisme  $\pi_{y}$ d\'efini sur les g\'en\'erateurs par :
$$\pi_{y}(f)= Y(l) \in V_{y} \subset Invd(T(L)) \  ,\  \pi_y (1)=id $$

$$\pi_y  :T (F) \rightarrow V_y (\Delta_F,\epsilon_F ) \subset Invd(T(L)) \subset
End(T(L)) \   .
$$
{\bf{D\'efinition.}}
\\
L'espace vectoriel des relations d'ordre $m$ de la big\`ebre $U_x$
est le sous espace vectoriel  $\ R_m (U_x ) \subset T_m (L )$ , noyau de l'application 
$$\pi_x :T_m (L) \rightarrow U_x \subset Invd(T(F)) \ \ . $$
On a de m\^eme les d\'efinitions similaires pour les espaces $R_n (V_y) \subset T_n (F)$
qui sont les noyaux des applications
$$\pi_y :T_n (F) \rightarrow V_y \subset Invd(T(L)) \ . $$ 

\subsection{Les espaces de relations d'ordre m sont des coid\'eaux.}

\vskip0.3cm
\begin{thm}
Les sous espaces vectoriels de relations $\ R_m (U_x )$,(d'ordre $m$) ,  respectivement
$\ R_n (V_y )$, (d'ordre $n$), sont des coid\'eaux des big\`ebres $T(L)(\Delta_L
,\epsilon_L )$ respectivement 
$T (F) (\Delta_F ,\epsilon_F )$ .\\  
Et les id\'eaux de relations pour les big\`ebres $U_x = \pi_x (T (L))$ et $V_y =\pi_y
(T (F ))$ sont des coid\'eaux. Ces id\'eaux sont engendr\'es par des sommes  de
coid\'eaux de dimensions finies .

\end{thm}
D\'emonstration.\\
Un \'el\'ement $w_m \in T_m (L) $ est une relation si et seulement si $\pi_x (w_m ) =0 $.
D'autre part, $\pi_x (w_m ) $ est un op\'erateur dans $U_x $ et en particulier un 
op\'erateur invariant \`a droite sur $T(F)$ ; il est nul si et seulement si :
$$ <\pi_x (w_m) , z> = 0 \ \ pour\ \ tout\ \ z\in T(F)\ . $$
 Le th\'eor\`eme de dualit\'e nous dit donc que  $w_m \in R_m (U_x ) $ si et seulement
si

$$<\pi_y (\tau (z) ) ,\tau (w_m )> =0 ,\forall z\in T (F) $$
o\`u $\tau $ d\'esigne l'anti-isomorphisme d'alg\`ebre tensorielle du renversement 
de l'ordre dans les tenseurs .\\
D'autre part $T_m (L) $ est une cog\`ebre de dimension finie, et les op\'erateurs 
invariants \`a droite $\pi_y (\tau (z))$ , laissent invariant cette cog\`ebre; ils
engendrent donc une sous-alg\`ebre de dimension finie dans le dual de $T_m (L) $; le
th\'eor\`eme  de dualit\'e exprime que $\tau (R_m (U_x )) $ est pr\'ecisement form\'e
par les
\'el\'ements de
$T_m (L)
$ qui s'annulent sur cette sous-alg\`ebre; $ \tau (R_m{U_x})   $ est donc un coid\'eal.
De plus parce que $\Delta_L :T(L) \rightarrow T(L)\otimes T(L) $ est un morphisme
d'alg\`ebres , l'on a:
$$\Delta _L \circ \tau = \tau \otimes \tau \circ \Delta_L \ ,$$
et
$$ \Delta_{L} (\tau (R_m )) \subset  \tau (R_m ) \otimes T_m (L) + T_m(L) \otimes \tau
(R_m ) \ ,    $$
et donc $ R_m (U_x ) $ est aussi un coid\'eal :
$$\Delta_L (R_m (U_x ))\subset R_m (U_x )\otimes T_m (L) +T_m (L) \otimes R_m
(U_x ) \ .$$
Il est d'autre part clair que la somme des id\'eaux engendr\'es par les espaces de
relations  d'ordre $m$ est l'id\'eal des relations de la big\`ebre
$U_{x}(\Delta_{L},\epsilon_L )$,
 et c'est bien un coid\'eal. Ce th\'eor\`eme valable pour toutes les big\`ebres
construites \`a partir de la donn\'ee d'une application lin\'eaire  $x:L\rightarrow
Invd(F) $ o\`u $L$ et $F$ sont deux cog\`ebres de dimensions finies est donc aussi
valable pour la big\`ebre $V_y (\Delta_F ,\epsilon_F )$ . \\ 
\\
Le th\'eor\`eme pr\'ec\'edent permet par relations d'orthogonalit\'es de d\'eterminer
l'ensemble des relations d'ordre $m$ pour chaque valeur de $m$, mais ne donne pas de
m\'ethode de les obtenir et surtout ne dit pas si l'id\'eal des relations est un id\'eal
engendr\'e par un 
nombre fini de relations .\\ Le th\'eor\`eme pr\'ec\'edent permet cependant de d\'ecrire
de fa\c{c}on pr\'ecise les conditions alg\'ebriques n\'ecessaires et suffisantes pour
l'obtention d'une ou plusieurs relations dans la big\`ebre 
$U_{x}(\Delta_{L},\epsilon_L )$,
\`a partir de la connaissance des relations dans l'alg\`ebre $Invd(F) $, des
relations dans l'alg\`ebre $K$ et de l'application $x:L\rightarrow Invd(F)$ . C'est ce
que nous allons d\'ecrire maintenant.\\
\\
Remarque.\\
\\
Soient $L$ et $F$ les cog\`ebres duales des alg\`ebres $K$ et $E$ et\\ $x:L \rightarrow
Invd(F)$; on rappelle que la donn\'ee de $x$ est exactement la donn\'ee de $$x_t :
L\rightarrow E \ .$$
Soit le produit des alg\`ebres $\otimes^{n}E$, $n \geq 0$ ; cette alg\`ebre est not\'ee :
$$\Pi (E) \  = \prod_{n\geq 0} \otimes^{n} E \ ,$$
avec $\otimes^{0} E =\mathbb{C} $, son dual est la cog\`ebre : 
$$T(F)$$
et si \`a $l\in L$  on associe l'\'el\'ement suivant, de l'alg\`ebre  $ \prod_{n\geq 0}
\otimes^{n} E \ $ ,
$$X_t (l)= (\otimes^{n}x_t \circ \Delta_{L}^{n-1}(l))_{n>0}$$
 avec pour  $n =0$  la valeur $1.\epsilon(l)$ .
Alors l'op\'erateur multiplication \`a gauche par $\ X_t (l)$ agit par transposition 
sur $T(F)$ et cette transposition co\"{\i}ncide avec l'op\'erateur $X(l)$.\\
\\
L'alg\`ebre $U_x $, engendr\'ee par les op\'erateurs $X(l)$ et l'identit\'e est donc
anti-isomorphe
\`a l'alg\`ebre engendr\'ee par les op\'erateurs $X_t (l)$ et la multiplication par
l'unit\'e de $\Pi(E)$.\\
\\
Le probl\`eme que l'on veux aborder est celui de l'obtention de conditions
alg\'ebriques n\'ecessaires et suffisantes pour la constuction des relations dans  $U_x
\subset Invd(T(F))$,
\`a partir de la connaissance des relations dans
$E$ ou de fa\c{c}on
\'equivalente dans $Invd(F)$, ainsi que des relations dans $K$, et de la donn\'ee de
l'application $x:L\rightarrow Invd(F)$ .\\
\\
Soit $w\in T_m (L)$, on consid\`ere dans $\Pi (K)$  la somme de tout les id\'eaux,
orthogonaux \`a $w$ . C'est le plus grand id\'eal orthogonal \`a $w$; et l'orthogonal
de cet id\'eal est la plus petite sous cog\`ebre $C_w $ de $T_m (L)$, contenant $w$.
De plus $\tau (C_w)$ est aussi une sous cog\`ebre de $T_m (L)$ qui contient $\tau (w)$ ,
et c'est la plus petite ; l'orthogonal dans $ \Pi (K)$ de $\tau (C_w)$ est donc le plus
grand id\'eal de $\Pi (K) $ orthogonal \`a $\tau (w)$ .\\
\\
Le th\'eor\`eme pr\'ec\'edent restant valable pour toutes sous cog\`ebres de $T
(L)$, parce que leurs images par $\tau $  restent de sous cog\`ebres et sont donc
invariantes sous l'action des op\'erateurs invariants
\`a droite sur $T(L)$ et en particulier sous l'action de la big\`ebre
$V_y(\Delta_F,\epsilon_F )$, l'on a :\\

\begin{thm}
Soient w une relation d'ordre $m$  et  $C_w$ la plus petite sous cog\`ebre
de $T_m(L)$ contenant $w$ ; $w\in R_m (U_x )\cap C_w
$ .\\ L'on a :\\ $ R_m (U_x )\cap C_w $ est un coid\'eal .
\end{thm} 
D\'emontrons maintenant le th\'eor\`eme suivant reliant les relations de $U_x$ aux
coid\'eaux de
$T (L)$ 
 contenus dans  $ker\ \pi_x^{1}  \cap ker\  \epsilon_L $ , $$ \pi_x^{1}
:T (L)\rightarrow Invd(F)\subset Invd(T(F))
$$

\begin{thm}
Soit $w\in R_m(U_x) $ , $w= \sum_{n\in I \subset(0...m)} w_n   $ , o\`u  $w_n \in
\otimes^n L$.
 $ R_m (U_x )\cap C_w $ est un coid\'eal de relations contenant $w$. 
L'orthogonal de ce coid\'eal est une sous alg\`ebre $B_w$ .Soit $B_w^{max}$ un
\'el\'ement maximal de l'ensemble des sous alg\`ebres contenues dans
$\Pi (K)$ , dont  l'orthogonal contienne $w$ et soit inclus dans $$ker\ \epsilon_L \cap  
ker\ (
\pi_x^{1} :T (L)\rightarrow Invd(F) )
$$
Alors l'orthogonal de  $B_w^{max}$  est un coid\'eal minimal de relations contenant $w$.
Ce coid\'eal est contenu dans le coid\'eal 
 $$ R_m (U_x )\cap C_w \ , $$   
o\`u $C_{w}$ est la plus petite sous cog\`ebre de $T (L)$ contenant $w$ .

\end{thm}
D\'emonstration.\\
L'orthogonal d'une telle alg\`ebre $B_w^{max}$ maximale, parceque c'est un coid\'eal de
$T (L)$ contenu dans $ker\ \epsilon_L \cap  ker\  (\pi_x^{1 }:T (L)\rightarrow
Invd(F) )
$, il est alors contenu dans l'id\'eal des relations de $U_x$.
En effet il suffit de d\'emontrer que pour tout $n \geq 2$ ,
$$\pi_{x}(u) (z) =0 \ \ pour z\in \otimes^{n}F $$
pour tout $u$ dans ce coid\'eal; l'on a parce que $U_x (\Delta_L ,\epsilon_L )$ agit sur
$T(F)$ :
$$\pi_{x}(u) (z) = \pi_x^{1} \otimes \pi_x^{1}\otimes ...\otimes \pi_x^{1}\circ
\Delta_{L}^{n-1}(u) (z)
\ ;$$  par hypoth\`ese $u$ appartient \`a un coid\'eal $C$ contenu dans\\ $ker\ 
\pi_{x}^{1}:T (L)
\rightarrow Invd(F)$ ; mais 
$$ \Delta_{L}^{n-1}(C) \subset C\otimes T (L)...\otimes T (L) + T (L)\otimes
C\otimes T(L)...\otimes T (L) +...+ T (L) ...\otimes T(L)\otimes C \ ,$$
ce qui d\'emontre que $ C$ est un coid\'eal de relations de $U_x $ parce que
$$\pi_{x}^{1} (C) = [0]\subset Invd(F) $$ par hypoth\`ese .\\
 D'autre part $B_w ^{max}$  \'etant un
\'el\'ement maximal, ce coid\'eal de relations est minimal et contient $w$ . \\
\\
Montrons que ce coid\'eal est
contenu dans
$ R_m (U_x )\cap C_w
$ . \\
Soit $I$ l'orthogonal  de $ C_w$ ; on rappelle que $I$ est le plus grand id\'eal de
$\Pi (K)$ orthogonal \`a $w$ ; on a : $I +B_{w}^{max}$ est une alg\`ebre dont
l'orthogonal contient $w$ et est contenu dans l'id\'eal des relations de $U_x $ ;
 parce que
$B_{w}^{max}$ est un \'el\'ement maximal, on a : 
$I +B_{w}^{max} = B_{w}^{max} $\\
et l'orthogonal de $ B_{w}^{max}$ est bien un coid\'eal contenu dans  $R_m \cap C_w $ 
.\\
\subsection{Premi\`ere construction de l'espace vectoriel des relations d'ordre m .}

Le th\'eor\`eme suivant d\'ecrit une m\'ethode pour obtenir par exemple l'espace
vectoriel des relations d'ordre $m$ de la big\`ebre $U_x(\Delta_L ,\epsilon_L )$ .\\
\\
Soit
le sous espace vectoriel $R_m^{1} $ de $T_m (L)$ d\'efini par le noyau de $\pi_x^{1}$ :
$$ R_m^{1} =T_m (L) \cap ker\  \pi_x^{1}\ \ ,\ \pi_x^{1} :T (L)\rightarrow
Invd(F)\subset Invd(T(F)) \ .$$ 
L'alg\`ebre  $\Pi (K)$ est le dual de $T (L)$ et l'orthogonal de $R_m^{1}$ est un
sous espace vectoriel de $\Pi (K)$ de la forme :
$$ W_m \oplus  \prod_{n>m} \otimes^{n} K        $$
o\`u  $W_m $ est un sous espace vectoriel de $T_m (K)$ .\\
Soit $ B(W_m)$ la plus petite sous alg\`ebre de $T_m (K)$ contenant $W_m$ et l'unit\'e
de l'alg\`ebre $T_m (K)$ ;\\ on a alors :
\begin{thm}
L'orthogonal dans $T_m (L)$ de la sous alg\`ebre $B(W_m)$ est $R_m(U_x)$ le coid\'eal
des relations d'ordre $m$ . 

\end{thm}
D\'emonstration :\\
L'orthogonal dans $T_m (L) $ de $ B(W_m)$ est un coid\'eal de $ T_m (L)$ contenu
dans $ker\epsilon_L \cap ker\  \pi_x^{1}$ ; c'est donc un coid\'eal de relations d'ordre
$m$ d'apr\'es le thm 4.3 ; d'autre part  on a vu, comme cons\'equence du th\'eor\`eme de
dualit\'e ,(thm 4.1), que l'espace des relations d'ordre $m$ , $R_m (U_x) $ est un
coid\'eal
\'evidemment contenu dans $ker\epsilon_L \cap ker\ \pi_x^{1}$ ; son orthogonal est une
alg\`ebre
$B_m$ qui contient
$W_m$et l'unit\'e de $T_m (K)$ ; donc $B(W_m) \subset B_m $ et par orthogonalit\'e on
obtient que l'espace  des relations d'ordre $m$ est contenu dans l'orthogonal de
$B(W_m)$, ce qui donne le th\'eor\`eme.\\
\\
Remarque.\\
Les th\'eor\`emes pr\'ec\'edents concernant les espaces de relations d'ordre $m$,
 $R_m (U_x) =T_m (L) \cap ker\ \pi_x \ \ ,\ \pi_x:T(L) \rightarrow U_x \subset Invd(T(F))
$ se g\'en\'eralisent pour toutes sous cog\`ebres de dimensions finies de $T (L)$ , 
en d\'efinissant l'espace des relations dans la cog\`ebre $C$ par :\\
 $R_C (U_x) =C \cap ker\ \pi_x  $ .
Ceci est d\^u aux faits que dans l'utilisation du th\'eor\`eme de dualit\'e, toute sous
cog\`ebre de $T (L) $ est invariante sous l'action des op\'erateurs $\pi_y (T(F))$,
qui sont des op\'erateurs invariants \`a droite sur $T(L)$, et que l'image par l'anti
-isomorphisme $\tau$ d'une cog\`ebre est une cog\`ebre .\\
\\
Pour obtenir une information plus pr\'ecise sur l'id\'eal des relations, il est
utile d'avoir une d\'escription plus concr\`ete de l'alg\`ebre donn\'ee par l'application
du th\'eor\`eme de dualit\'e dans le th\'eor\`eme 4.1 .

\subsection{Construction compl\'ementaire de l'espace vectoriel des relations d'ordre m.}

Le th\'eor\`eme de dualit\'e nous a permis de montrer que l'espace vectoriel des
relations $R_m(U_x)$ est un coid\'eal de la cog\`ebre $T_m(L)$ , dont le dual est
$T_m(K)$.  Soit $B_m \subset T_m(K) $ l'orthogonal de  $R_m(U_x)$; c'est une sous
alg\`ebre de $T_m(K)$.  Le th\'eor\`eme pr\'ec\'edent a donn\'e une caract\'erisation
de cette alg\`ebre:\\
$B_m$ est l'alg\`ebre engendr\'ee par l'espace vectoriel $W_m$ et l'unit\'e de $T_m (K)$
.\\
$W_m$ est l'orthogonal de l'espace vectoriel $R_m^{1}$ :
 $$R_m^{1}= T_m(L) \cap  ker (\pi_x^{1} :T (L)\rightarrow
Invd(F) ) \ .$$\\
Nous allons montrer que $B_m $ co\"{\i}ncide avec l'alg\`ebre $B^m$ engendr\'ee de la
fa\c{c}on suivante dans $T_m (K)$:\\
\\
pour tout $f\in F$ soit $Y_{t,op}^m(f)$ d\'efini par :\\
$$Y_{t,op}^m (f) = \sum_{n\in (0..m)}\otimes^n y_{t} \circ \Delta_{F,op}^{n-1}(f) \in T_m
(K)
$$ o\`u $\Delta_{F,op}$ d\'esigne le coproduit oppos\'e sur $F$ ,\\
et o\`u  $y:F\rightarrow Invd(L) $ d\'efinit $y_t :F \rightarrow K$ par la
co\"{\i}ncidence de $y(f)$ avec la transposition de l'op\'erateur multiplication \`a
gauche par $y_t (f)$ ; pour $n=0$ , $Y_{t,op}(f) = \epsilon_F (f)$.\\
\\ Soit
$B^m$ l'alg\`ebre engendr\'ee dans $T_m (K)$ par l'espace vectoriel
$Y_{t,op}^m (F)$ et l'identit\'e de $T_m (K)$ ;\\ on a le th\'eor\`eme suivant:

\begin{thm}
 L'orthogonal de $B^m$ dans $T_m(L)$ est l'espace des relations d'ordre $m$ :
$R_m(U_x)$ ;  et l'on a donc  $$B_m =B^m \ \ ,\ \ \forall m  \ .$$

\end{thm}
D\'emonstration.\\
Dans le th\'eor\`eme 4.1 on a vu que $w\in T_m(L) $ est une relation de $U_x$
si et seulement si :
$$<\pi_y(T(F)),\tau (w)>= 0$$
si et seulement si 
$$\epsilon_L \circ \pi_y(T(F)) (\tau (w)) = 0 \ .$$ 
En particulier \'evaluons cette expression sur $f\in F \subset T(F)$ ; par la
lin\'earit\'e de l'op\'erateur $\pi_y (f)$ sur $T(L)$ , l'on a :

$$\epsilon_L \circ \pi_y(f) (\tau (w)) =\sum_{n\in (0,m)} <\otimes^n 
y \circ \Delta_F^{n-1} (f),\tau (w_n), > =\sum_{n\in (0,m)}< \otimes^n y \circ \tau
\circ 
\Delta_{F}^{n-1}(f) ,w_m > \ ,$$
ce qui est \'egal \`a
$$ \sum_{n\in (0,m)}< \otimes^n y  \circ 
\Delta_{F,op}^{n-1}(f) ,w_m > 
=<\pi_{y,op}^m(f), w>.$$ o\`u $  \pi_{y,op}^m(f)$ est l'op\'erateur invariant \`a droite
sur  $T_m (L)$ correspondant \`a la multiplication \`a gauche par $Y_{t,op}^m (f)$ dans
$T_m (K) $.\\
\\
Consid\'erons plus explicitement l'op\'erateur 
$\pi_y (f) \circ \tau $ :
$$\pi_y (f) \circ \tau = \epsilon_L \circ \pi_y(f) \otimes id \circ \tau \otimes \tau
\circ \Delta_L 
$$ 
o\`u l'on a utilis\'e le fait que $\Delta_L$ et $\tau $ sont respectivement un
morphisme d'alg\`ebres et un anti-homomorphisme d'alg\`ebre.\\
On obtient donc :
$$\pi_y (f) \circ \tau = \tau \circ \pi_{y,op}(f) . $$ 
Donc $ w\in T_m (L) $ est une relation d'ordre $m$ de la big\`ebre $U_x$
si et seulement si pour toute suite $(f_1,f_2,...f_q ) $ l'on a :
$$ \epsilon_L \circ \pi_{y}(f_{q})\circ ...\circ  \pi_{y}(f_{2})\circ 
\pi_{y}(f_{1})\circ \tau (w) =0 \ et \ \epsilon_L \circ \pi_y (1)(\tau(w))
=0=\epsilon_L (w) .   
$$ Par l'it\'eration de la formule pr\'ec\'edente on obtient donc que  :
$w\in R_m (U_x)$ si et seulement si pour toute suite $(f_1,f_2,...f_q ) $  l'on a :
$$ \epsilon_L \circ \tau \circ \pi_{y,op}(f_{q})\circ ...\circ  \pi_{y,op}(f_{2})\circ 
\pi_{y,op}(f_{1}) (w) =0  \ \ et\ \ \epsilon_L (w) =0\ .  $$
L'on a aussi $ \epsilon_L \circ \tau = \epsilon_L $, et utilisant l'anti-isomorphisme 
entre les op\'erateurs invariants \`a droite sur $T_m(L)$ et l'alg\`ebre des
op\'erateurs multiplications \`a gauche sur $T_m (K)$ l'on obtient :\\
\\ 
 $w$ est une relation d'ordre $m$ de $U_x$ si et seulement si $w$ est orthogonal
\`a l'alg\`ebre $B^m$ engendr\'ee dans $T_m(K)$ par les op\'erateurs $Y_{t,op}^m (F)$
et l'identit\'e de $T_m (K)$.\\
\\
La compl\'ementarit\'e de ces deux constructions de $R_m(U_x)$ provient de ce que les
deux alg\`ebres identiques qui d\'efinissent l'orthogonal de cet espace de relations,
$B_m$ et $B^m$ sont engendr\'ees par des sous espaces vectoriels diff\'erents.
\vskip7cm
{\bf{R\'ef\'erences}}
\vskip0.4cm 
1) E. Mourre : Remarques sur une m\'ethode de r\'ealisations de big\`ebres ,

 et
alg\`ebres de Hopf associ\'ees \`a certaines r\'ealisations ,
arXiv:math-ph/0312049.\\
\vskip0.3cm 
2) E. Abe : Hopf algebras , Cambridge University Press , 1980 .\\
 \vskip0.3cm 
3) C. Kassel , M. Rosso , V. Turaev :  Quantum groups and

     knots invariants , Panoramas et Synth\`eses , 97 , (N. 7) , S.M.F.
\vskip0.4cm 
4)  V.G. Drinfeld : Quantum groups , proc.int.conf.math. , Berkley,California, 

1986
.\\
 \vskip0.3cm  
5)  S. Majid : Foundations of quantum groups theory ,

 Cambridge
University Press, 1995.\\ 
\vskip0.3cm 
6) A. Connes , H. Moscovici : C.M.P. 198 1998 .
\vskip0.3cm
7) A. Connes , D. Kreimer  :Ann.Henri Poincar\'e . (V 3) 2002 .
\vskip0.3cm
8) A.Connes , M. Marcolli  : ArXiv: hep-th 0411114 v1 2004 .

\end{document}